\def\[{[\![}
\def\]{]\!]}
\def\l{\langle}
\def\r{\rangle}
\def\Z{{\bf Z}}
\def\C{{\bf C}}
\def\Ch{{\C[[h]]}}
\def\N{{\bf N}}

\def\s{\smallskip}

\def\n{\noindent}
\def\Uq{U_q[osp(2n+1|2m)]}

\def\q{{\bar q}}

\def\L{{\bar L}}

\def\a{\alpha}

\def\t{\theta}

\def\tL{{\tilde L}}
\def\0{{\dot 0}}

\baselineskip=16pt

\font\twelvebf=cmbx12

\vskip 20mm

\noindent {\twelvebf Deformed Clifford Cl$_{\bf q}$(n$|$m)
and orthosymplectic U$_{\bf q}$[osp(2n+1$|$2m)] superalgebras
and their root of unity representations }

\vskip 32pt
H.-D. Doebner$^1$, T. D. Palev$^2$\footnote*{Permanent address: 
Institute for Nuclear
Research and Nuclear Energy, 1784 Sofia, Bulgaria; E-mail:
tpalev@inrne.bas.bg, stoilova@inrne.bas.bg}, N. I. Stoilova$^1$*

$^1$Department of Physics, Technical University
of Clausthal

Leibnizstrasse 10, D-38678 Clausthal-Zellerfeld, Germany

\smallskip

$^2$Abdus Salam International Centre for Theoretical Physics, 
PO Box 586, 34100 Trieste, Italy

\smallskip
E-mail: asi@pt.tu-clausthal.de, palevt@ictp.trieste.it, 
ptns@pt.tu-clausthal.de

\vskip 2cm

\n{\bf Abstract.} It is shown that the Clifford superalgebra $Cl(n|m)$ 
generated by
$m$ pairs of Bose operators (odd elements) anticommuting with $n$
pairs of Fermi operators (even elements) can be deformed to
$Cl_q(n|m)$ such that the latter is a homomorphic image
of the quantum superalgebra $U_q[osp(2n+1|2m)]$. The Fock space
$F(n|m)$ of $Cl_q(n|m)$ is constructed. At $q$ being a
root of unity ($q=\exp ({i\pi l/k})$) $q-$bosons (and 
$q-$fermions) are operators acting in a finite-dimensional
subspace $F_{l/k}(n|m)$ of $ F(n|m)$. Each $F_{l/k}(n|m)$
is turned through the above mentioned homomorphism
into an irreducible (root of unity) $U_q[osp(2n+1|2m)]$
module.
For $q$ being a primitive root of unity ($l=1$) the corresponding 
representation is 
unitary. The module
$F_{1/k}(n|m)$ is decomposed into a direct sum of irreducible
$U_q[sl(m|n)]$ submodules. The matrix elements of all Cartan-Weyl
elements of $U_q[sl(m|n)]$ are given  within each such
submodule.

\vfill\eject

\leftskip 0pt

\n
{\bf 1. Introduction}

\bigskip
\n
{\bf 1.1 Some definitions}

\smallskip\n
In this subsection we give definitions which we need to present 
the motivations and the aims of the paper.

Let $(a_{ij})$ 
[1], be an $(m+n) \times (m+n)$ symmetric Cartan
matrix
with entries:
$$ (a_{ij})=(-1)^{\l j\r}
\delta_{i+1,j}+(-1)^{\l i \r}
\delta_{i,j+1} -[(-1)^{\l j+1
\r} + (-1)^{\l j \r}]\delta_{ij}+
\delta_{i,m+n}\delta_{j,m+n}, \eqno (1)
$$
where
$$
\langle i \rangle = \cases
{\bar 1, & for $i\le m$ \cr
\bar 0, &  for $i> m$ \cr}  \eqno(2)
$$
and  $\Z_2=\{\bar{0},\bar{1}\}$ is the ring of all integers modulo 2.

\bigskip
\n {\bf Definition 1}  
[1]. $U_q[osp(2n+1|2m)]$ is a Hopf algebra,
which is a topologically free $\Ch$ module, $q=e^h$, (complete in the
$h-$adic topology), 
with Chevalley generators $h_i, e_i,f_i$, $
i=1,\ldots ,n+m=N,$ subject to the
\n
1. Cartan-Kac relations:
$$
\eqalignno{
& \; [h_i,h_j]=0, & (3a) \cr
& \; [h_i,e_j]=a_{ij}e_j,& (3b) \cr
& \; [h_i,f_j]=-a_{ij}f_j, & (3c) \cr
& \;  \[e_i,f_j\]=\delta_{ij}{{k_i-{\bar k}_i}\over{q-{\bar q}}},~~~
 q=e^h,~{\bar q}=q^{-1},~k_i=q^{h_i},~ {\bar k}_i = k_i^{-1}=q^{-h_i},  &
(3d) \cr
}
$$
2. $e-$Serre relations
$$
\eqalignno{
& \[e_i,e_j\]=0, \quad \vert i-j \vert \neq 1, \quad e_m^2=0,& (4a) \cr
& [e_i,[e_i,e_{i \pm 1}]_{\bar q}]_q \equiv
  [e_i,[e_i,e_{i \pm 1}]_q]_{\bar q}=0, \quad
  i\neq m,\;\; i\neq m+n, & (4b) \cr
& \{[e_{m},e_{m-1}]_{q},[e_m,e_{m+1}]_{\bar q}\}
\equiv \{[e_{m},e_{m-1}]_{\q},[e_m,e_{m+1}]_{q}\}=0, & (4c)\cr
& [e_N,\[e_N,[e_N,e_{N-1}]_{\bar q}\]]_q \equiv
[e_N,\[e_N,[e_N,e_{N-1}]_q\]]_{\bar q}=0, & (4d)\cr
}
$$

\noindent
3. $f-$Serre relations $\hfill (5)$

\n
obtained from (4) by
replacing everywhere $e_i$ with $f_i$.

\n
The grading
is induced from:
$${\rm deg}(h_j)=\bar 0, \;\;\forall j, \;\;
{\rm deg}(e_m)={\rm deg}(f_m)=\bar 1 ,
\quad {\rm deg}(e_i)={\rm deg}(f_i)=\bar 0
\quad {\rm for}\quad i\neq m.
$$
Above and throughout
$\C[[h]]$ is the ring of all
complex formal power series in $h$;

\n
$
[a,b]=ab-ba,\;\; \{a,b\}=ab+ba,
\;\;\[a,b\]=ab-(-1)^{deg(a)deg(b)}ba; \hfill (6)
$

\n
$
[a,b]_x=ab-xba,\;\; \{a,b\}_x=ab+xba,
\;\;\[a,b\]_x=ab-(-1)^{deg(a)deg(b)}xba. \hfill (7)
$

We do not write the other Hopf superalgebra maps
(comultiplication, co-unit and antipode) since we will not use
them. Certainly they  are also a part of the definition.

An alternative description of the algebra under consideration
in terms of the so called deformed Green generators was given in Ref.
[2]:

\n
{\bf Definition 2.} $\Uq$ is a topologically free $\Ch$ module
($q=e^h$)
and an associative unital 
superalgebra with generators $H_i,\;
a_i^\pm,\;\;i=1,\ldots,m+n=N$ and relations:

$$
\eqalignno{
& [H_i,H_j]=0,
& (8a) \cr
& [H_i,a_j^\pm]=\pm \delta_{ij}(-1)^{\l i \r} a_j^\pm ,
& (8b)\cr
& \[a_i^-,a_i^+\]=-2{L_i-{\bar L}_i\over q-{\bar q}},
~~ L_i=q^{H_i}, \quad\L_i=q^{-H_i},
& (8c) \cr
& [\[a_{N-1}^\xi,a_N^\xi\],a_{N}^\xi]_{q^{-(-1)^{\l N\r}}}=0,
& (8d) \cr
& \[\[a_i^{\eta},a_{i+ \xi}^{-\eta}\],
a_j^{\eta}\]_{q^{-\xi (-1)^{\l i \r}\delta_{ij}}}
=2(\eta)^{\l j \r}\delta _{j,i + \xi}L_j^{-\xi \eta}a_i^{\eta},&(8e)
\cr
}
$$
where
$$
{\rm deg}(a_i^\pm)\equiv \langle i \rangle . \eqno(9)
$$
The expressions of the deformed Green generators via the
Chevalley generators read
($i=1,\ldots,N-1$):
$$
\eqalignno{
&  a_i^-=(-1)^{(m-i)\l i \r}\sqrt{2}[e_i,[e_{i+1},
  [\ldots,[e_{N-2},[e_{N-1},e_N]_{q_{N-1}}]_{q_{N-2}}
  \ldots ]_{q_{i+2}}]_{q_{i+1}}]_{q_{i}},\quad
   a_N^-=\sqrt{2}e_N, & (10a)   \cr
&  a_i^+=(-1)^{N-i+1}\sqrt{2}
   [[[\ldots [f_N,f_{N-1}]_{{\bar q}_{N-1}},f_{N-2}]_{{\bar q}_{N-2}}
   \ldots]_{{\bar q}_{i+2}},f_{i+1}]_{{\bar q}_{i+1}},f_{i}]_{{\bar
q}_{i}}
   , \quad a_N^+=-\sqrt{2}f_N,  &(10b) \cr
&  H_i=h_i+h_{i+1}+\ldots +h_N \;({\rm including} \;i=N),& (10c)  \cr
}
$$
with
$$
q_i=q^{(-1)^{\l i+1 \r}},\quad  {\it i.e.}~~~
q_i={\bar q},\;i<m, \quad q_i=q,\;i\ge m. \eqno(10d)
$$
For the inverse relations one finds:
$$
\eqalignno{
& e_i={1\over 2}\L_{i+1}\[a_i^-,a_{i+1}^+\],
\quad e_N={1\over \sqrt{2}}a_N^-, & (10e) \cr
& f_i=-{1\over 2}(-1)^{\l i+1 \r}\[a_{i+1}^-,a_{i}^+\]L_{i+1}
     ={1\over 2}\[a_i^+,a_{i+1}^-\]L_{i+1} ,
\quad f_N=-{1\over \sqrt{2}}a_N^+, & (10f) \cr
& h_i=H_i-H_{i+1}, \quad H_N=h_N. & (10g)\cr
}
$$
In the above definitions $h$ is an abstract indeterminate. We
shall be dealing however with representations of $\Uq$ where $h$
is a complex number. The same holds also for the Clifford superalgebra
defined in Section 2 and for $U_q[sl(m|n)]$. 

\bigskip\n 
{\bf 1.2 Motivations and Aims}

\smallskip\n
The motivation for giving a definition of $\Uq$ in terms of 
the deformed Green generators (definition 2)
stems from the observation that the nondeformed
Green generators $\hat{a}_i^\pm,\;\;i=1,\ldots,m+n,$ of
$osp(2n+1|2m)$, contrary to the Chevalley generators, have a
physical meaning [3]: $\hat{a}_1^\pm,\ldots ,
\hat{a}_m^\pm$ (resp. $\hat{a}_{m+1}^\pm,\ldots ,
\hat{a}_{m+n}^\pm$ ) are para-Bose (pB) (resp. para-Fermi (pF))
operators [4]. In one particular representation of $osp(2n+1|2m)$
the para-Bose (resp. para-Fermi) operators become usual Bose
(resp. Fermi) operators. 
Moreover the grading of these operators is
opposite to the canonical one: the Bose operators ($\equiv$
the bosons) are odd (fermionic) variables, 
the Fermi
operators ($\equiv$ the fermions) are even (bosonic) variables, 
and the bosons anticommute with the fermions. Let us
underline that the above somewhat unusual properties are not an
input. On the contrary, they are an output from the
representation theory of $osp(2n+1|2m)$.

The above result can be viewed also in the following way. Let
$Cl(n|m)$ be the associative superalgebra which is generated by
(i) $m$ pairs of Bose operators ${\hat c}_1^\pm,\ldots,{\hat
c}_m^\pm$  (postulated to be odd variables), (ii) $n$ pairs of
Fermi operators ${\hat c}_{m+1}^\pm, \ldots,{\hat c}_{m+n}^\pm$ (even
elements) (iii) and the assumption that the bosons anticommute
with the fermions. Then the linear map $\varphi$ defined on the
generators as
$$
\varphi({\hat a}_i^\pm)= {\hat c}_i^\pm, ~~i=1,\ldots,n+m \eqno(11)
$$
and extended by associativity is a homomorphism of the universal
enveloping algebra (UEA) $U[osp(2n+1|2m)]$ onto $Cl(n|m)$.

Independently from [3], bosons and fermions with the above
properties (in the case $n=m=\infty$) were used in [5] for
studying of irreducible representations of some
infinite-dimensional superalgebras (including
$osp(\infty|\infty)$ and $gl(\infty|\infty)$), relevant for
construction of supersymmetric generalizations of certain
hierarchies of soliton-like evolution equations.
Accepting the 
terminology of this paper we refer to $Cl(n|m)$ as to Clifford
superalgebra.

There are different approaches for extending  the results of
single mode deformed bosons $b^\pm $ 
 [6-8] or fermions $f^\pm $ 
[9] to the multi-mode case $b_i^\pm, f_i^\pm, i=1,2,\ldots $.
One way is to preserve the main features of the
nondeformed case, postulating that different modes of $q-$bosons
commute and different modes of $q-$fermions anticommute (see [9]
for the precise definitions). In other approaches the relations
between different modes are not postulated. They are derived on
the ground of other assumptions, yielding as a rule that
different modes of operators do not commute, but "$q-commute$''
[10-13].

Our approach for (simultaneous) deformation of bosons and
fermions is based on the requirement the Eq. (11) to be
preserved also in the quantum case. More precisely, we require
that the relations between the generators of the deformed
Clifford superalgebra $Cl_q(n|m)$, namely the relations between
$m$ pairs of deformed bosons ($q-$bosons)
$c_1^\pm,\ldots,c_m^\pm$ and  $n$ pairs of deformed fermions
$q-$fermions
$c_{m+1}^\pm,\ldots,c_{m+n}^\pm$ are determined in
such a way that the $q-$analogue of the above map (11) is a
homomorphism of $\Uq$ onto $Cl_q(n|m)$.
As a result the single mode
$q$-bosons (resp. $q$-fermions) 
coincide (up to an overall
multiple) with the known $q$-bosons
[6-8,14], but different modes of them ``$q$-commute'' (resp.
``$q$-anticommute''). 
The deformed bosons and fermions
mutually $q$-anticommute. The $q-$bosons 
$c_1^\pm,\ldots , c_m^\pm$ coincide with those introduced in [15]
in connection with the quantum superalgebra $U_q[osp(1|2n)]$.

One of the consequences of the present approach stems from the
observation that any (irreducible) $Cl_q(n|m)$ module $V$ and in
particular its Fock space $F(n|m)$ can be immediately
turned into an (irreducible) representation space of $\Uq$ and
also of representation space of any other subalgebra of it, in
particular of $U_q[sl(m|n)]$. Based on this we shall construct
root of unity irreducible representations of $\Uq$ and of
$U_q[sl(m|n)]$ in (appropriate subspaces/factorspaces of) 
$F(n|m)$.

In Section 2 we introduce the deformed Clifford algebra
$Cl_q(n|m)$ and study its Fock space $ F(n|m)$. The
transformations of the basis under the action of the $Cl_q(n|m)$
generators are written down (Eqs. (15)). Next we show that when
$q=\exp({i\pi l/k})$ with $l,\, k$ - integers, i.e. when $q$ is a
root of unity, $F(n|m)$ contains an infinite-dimensional
invariant (and depending on $q$) subspace $I_{l/k}$, so that the
factor space $F_{l/k}(n|m)=F(n|m)/I_{l/k}$, is a
finite-dimensional irreducible $Cl_q(n|m)$ module. We show that
the representation of $Cl_q(n|m)$ in $F_{l/k}(n|m)$ is unitary
(with respect to a natural antiinvolution) only if $q$ is a
primitive root of unity ($l=1$). The transformations of the basis
of $F_{1/k}(n|m)$ under the action of the deformed Bose and Fermi
operators are written down explicitly (Eqs. (46)).

Using the homomorphism $\varphi$ of $U_q[osp(2n+1|2m)]$ onto
$Cl_q(n|m)$ (see Proposition 6), in Section 3 we turn each Fock
spaces $F_{l/k}(n|m)$
into an irreducible $U_q[osp(2n+1|2m)]$ module,
thus finding a class of representations of this quantum superalgebra.
Next we consider $U_q[sl(m|n)]$ as a subalgebra of
$U_q[osp(2n+1|2m)]$. As one can expect, the Fock spaces carry
reducible representations of $U_q[sl(m|n)]$. Each Fock space
$F_{1/k}(n|m)$ is decomposed into a direct sum of $U_q[sl(m|n)]$
irreducible subspaces. We write down the transformations of these
irreducible $U_q[sl(m|n)]$ modules not only with respect to the
Chevalley generators but under the action of all Cartan-Weyl
elements of $U_q[sl(m|n)]$.

\bigskip\bigskip
\n
{\bf 2. A q-deformed Clifford superalgebra Cl$_q$(n$|$m) and its Fock
representations}

\bigskip\n
In this section we define a $q$-deformed Clifford superalgebra
$Cl_q(n|m)$, generated by a set of operators interpreted as
deformed Bose and Fermi operators. We study a class of Fock
representations of $Cl_q(n|m)$, which are finite-dimensional in
the case when $q$ a root of unity. In this way parallel to the
$q$-fermions also the $q-$bosons are finite-dimensional.

The restriction $q$ to be a root of unity is not artificial. It
is a consequence of basic physical requirements (e.g. the physical
observables have to be Hermitian operators). Mathematically this
leads to representations, which are unitary with respect to a
natural antilinear anti-involution (see Eq. (21)). It turns out
that only the representations corresponding to $q$ being a
primitive root of unity are unitary.

The results of the present section will be used in the next one
in order to define root of unity representations of the superalgebra
$\Uq$ and its subalgebra $U_q[sl(m|n)]$.

Define a set of  $\Z_2$-graded
operators $c^\pm_i, ~~ N_i$
$$
 deg(c_i^\pm)=\l i \r, ~~ deg(N_i)=\bar{0},
  ~~i=1,\ldots , m+n
$$
with relations:
$$
\eqalignno{
& [N_i, N_j]=0,
& (12a) \cr
& [N_{i},  c^{\pm}_j] = \pm \delta_{ij} c^\pm_j ,
& (12b) \cr
& c_i^-c_i^+ +(-1)^{\langle i \rangle}q^{\pm 1}c_i^+c_i^-
  ={2\over {q^{1/2}+q^{-1/2}}}q^{\pm (-1)^{\l i \r}N_{i}}, & (12c)
\cr
& c_i^\xi c_j^\eta =-(-1)^{\l j \r}q^{\xi \eta}
  c_j^\eta c_i^\xi , \quad {\rm for\ all\ }i<j,~~
\xi, \eta =\pm ~~or~ \pm 1,
 & (12d) \cr
& (c_i^{\pm})^2=0, ~~ i=m+1,\ldots ,m+n. & (12e) \cr
}
$$
Note that Eqs. (12c) are equivalent to
$$
\eqalignno{
& c_i^+c_i^-= c_q{q^{N_i}-\bar{q}^{N_i}
\over{q-\bar{q}}},  ~~c_q={2\over{q^{1/2}+\bar{q}^{1/2}}}, & (12f) \cr
& c_i^-c_i^+= c_q{q^{1-(-1)^{\l i \r}N_i}-\bar{q}^{1-(-1)^{\l i \r}N_i}
\over{q-\bar{q}}}.
& (12g)\cr
}
$$
It is obvious that in the limit $q \rightarrow 1$ the operators
$c_1^\pm,\ldots,c_m^\pm$ (resp. $c_{m+1}^\pm,\ldots,c_{m+n}^\pm$)
reduce to Bose (resp. Fermi) creation and annihilation operators,
which mutually anticommute, and $N_i$ are the corresponding
number operators.

\bigskip\noindent 
{\bf Definition 3.} The deformed Clifford algebra $Cl_q(n|m)$
is a topologically free $\Ch$ module and an associative
unital superalgebra
with generators $c_i^{\pm},$ $ N_{i}$ subject
to the relations (12).

Clearly, $Cl_q(n|m)$ is a deformation of the Clifford
superalgebra $Cl(n|m)$.

\bigskip
We proceed with the construction of the Fock space 
$F(n|m)$ of $Cl_q(n|m)$. The vacuum vector $|0\r $ is defined in a
natural way:
$$
c_i^-|0\r =0, \quad
N_{i}|0\r =|0\r, \quad
i=1,\ldots,m+n.
\eqno(13)$$
As a basis take
$$
(c_1^+)^{r_1}(c_2^+)^{r_2}\ldots (c_{m+n}^+)^{r_{m+n}}|0\r=
|r_1,r_2, \ldots , r_{m+n} \r ,~
r_i\in {\bf Z_+}, ~i=1,\ldots ,m; ~r_i\in \{ 0,1 \},~i=m+1,
\ldots ,m+n, \eqno(14)
$$
where $\Z_+$ are all non-negative integers.

\n
{\bf Proposition 1.} The transformation of the basis (14) under
the action of the deformed Bose and Fermi operators reads:
$$
\eqalignno{
&N_{i}|r_1,\ldots , r_{m+n} \r =
r_i|r_1,\ldots , r_{m+n} \r, & (15a)\cr
&&\cr
&c_i^+|r_1,\ldots , r_{m+n} \r =
(-1)^{(1-\l i\r )(r_1+\ldots +r_{i-1})}(1-(1-\l i \r )r_i)
\bar{q}^{r_1+\ldots +r_{i-1}}
|\ldots , r_{i-1}, r_i+1, r_{i+1}, \ldots
\r ,& (15b)\cr
&&\cr
&c_i^-|r_1,\ldots , r_{m+n}\r =
(-1)^{(1-\l i \r )(r_1+\ldots +r_{i-1})}
[r_i] c_q q^{r_1+\ldots +r_{i-1}}
|\ldots , r_{i-1}, r_i-1, r_{i+1}, \ldots, \r,
& (15c)\cr
}
$$
where
$$
[x]={q^x-q^{-x}\over{q-\bar{q}}}.
\eqno(16)
$$

\bigskip
The proof follows from
the defining relations (12) and identities like:
$$
\eqalignno{
& 1+q^2+q^4+\ldots +q^{2n}+\bar{q}^{2}+\bar{q}^{4}+\ldots +\bar{q}^{2n}=[2n+1],
 & (17)\cr
&&\cr
& q+q^3+\ldots +q^{2n+1}+\bar{q}+\bar{q}^{3}+\ldots +\bar{q}^{2n+1}=[2n+2],
& (18)\cr
&&\cr
&c_i^-(c_i^+)^n=\bar{q}^n(c_i^+)^nc_i^-+c_q[n](c_i^+)^{n-1}q^{N_i},
~~
i=1,\ldots ,m.
& (19)\cr
}
$$
\hfill $[]$

On the ground of the Fock space $F(n|m)$ we proceed to
determine a class of unitary representations [16] of $Cl_q(n|m)$ with
respect to the antilinear anti-involution $\omega$ defined by
$$
\omega(c_i^\pm)= c_i^\mp,~~~\omega(N_i)=N_i. \eqno(20)
$$
If (~,~) is a nondegenerate Hermitian form then (20) is
equivalent to
$$
\eqalignno{
&( c_i^+)^{\dag} =c_i^-, \quad ( N_{i})^{\dag} =N_{i},
\quad i=1,\ldots, m+n,
&(21) \cr
}
$$
with $x^{\dag}$ being the conjugate to $x$ operator with respect to
$(~,~)$ and the corresponding representation is said to be
contravariant. If (~,~) is a scalar product, the representation is
unitary.

The unitarity condition (21) stems from physical considerations.
For example in  case of a harmonic oscillator it is convenient to
replace the position operators $q_k$ and the momentum $p_k$ via
creation and annihilation operators $c_k^\pm$: $q_k=(c_k^+ +
c_k^-)/{\sqrt 2}$, $p_k=i(c_k^+ - c_k^-)/{\sqrt 2}$. Then the
necessity $q_k$ and $p_k$ to be Hermitian operators leads to Eq.
(21). For the same reason the unitarity condition has to hold in
any second quantized picture, in the filling number
representation, etc. Let us add that the condition (21) is
statistically independent, it has to hold for any statistics
compatible with the principles of quantum theory.

Define a Hermitian form (~,~) on $F(n|m)$  requiring
$$
( |0\r, |0\r ) = \l 0|0 \r =1 \eqno (22)
$$
and postulating that (21) should be satisfied, namely
$$
( c_{i}^{\pm}v,w) =( v, c_i^{\mp}w ), ~~( N_{i}v,w) =( v, N_i w ), ~~
v,w\in F(n|m).
 \eqno (23)
$$

\noindent
{\bf Proposition 2.} The above definition is selfconsistent
only if $q$ is a pure phase, more precisely if
$$
q=\exp ({i\phi}), ~~~0\le \phi <2\pi ,~~~\phi\ne \pi. \eqno(24)
$$
{\it Proof.} Only for simplicity we assume that $m>1$. From (12),
(13), (22) and (23) one derives uniquely that
$$
(|1,{\dot 0}\r,|1,{\dot 0}\r)=(|0,1,{\dot 0}\r,|0,1,{\dot 0}\r)=c_q,
\eqno(25)
$$
where ${\dot 0}$ indicates that all remaining entries are zero,
$|1,{\dot 0}\r=|1,0,\ldots,0\r$.
$c_q$ has to be a nonzero real number. Compute now
$
(|1,1,{\dot 0}\r,|1,1,{\dot 0}\r).
$
$$
\eqalign{
&(|1,1,{\dot 0}\r,|1,1,{\dot 0}\r)=(c_1^+c_2^+ |0\r,c_1^+ c_2^+|0
\r)=
(c_2^+ |0\r,c_1^- c_1^+ c_2^+|0 \r)\cr
&=(c_2^+ |0\r, c_q\ {q^{1+N_1}-\bar{q}^{1+N_1}
\over{q-\bar{q}}}\ c_2^+|0 \r)=c_q(|0,1,{\dot 0}\r,|0,1,{\dot
0}\r)=c_q^2.\cr
}\eqno(26)
$$
On the other hand the 
LHS of (26) can be evaluated in a different
way:
$$
\eqalign{
& (|1,1,{\dot 0}\r,|1,1,{\dot 0}\r)=
(c_1^+c_2^+ |0\r,c_1^+ c_2^+|0 \r)=(q c_2^+c_1^+ |0\r,q c_2^+ c_1^+|0
\r)
=|q|^2 (c_1^+ |0\r,c_2^- c_2^+ c_1^+|0 \r)\cr
&=|q|^2(c_1^+ |0\r, c_q{q^{1+N_2}-\bar{q}^{1+N_2}
\over{q-\bar{q}}}c_1^+|0 \r)=|q|^2 c_q(|1,{\dot 0}\r,|1,{\dot 0}\r)
=|q|^2 c_q^2.\cr
}\eqno(27)
$$
Eqs. (26) and (27) are compatible only for $|q|=1$. From the
latter and (12f) it follows that $c_q=\Big(cos(\phi/2)\Big)^{-1}$
and therefore $\phi \ne \pi$. Hence (24) holds.

\hfill $[]$

\bigskip
Any two vectors $|r_1,\ldots , r_{m+n} \r $
~~and ~~~$|r_1^{\prime},\ldots , r_{m+n}^{\prime}
\r $ ~~~with ~~~
$(r_1,\ldots , r_{m+n}) $ $\neq $
$(r_1^{\prime},\ldots , r_{m+n}^{\prime})$ are orthogonal and
$$
(|r_1,\ldots , r_{m+n}\r ,
|r_1,\ldots , r_{m+n} \r )=c_q^{r_1+\ldots +r_{m+n}}
[r_1]![r_2]!\ldots [r_m]!, \eqno(28)
$$
where $[r]!=[r][r-1]\ldots [1]$.
Hence $(~,~)$ is a nondegenerate Hermitian form and therefore Eq. (21)
is a consequence of (20).

\bigskip
\n
{\bf Proposition 3.} The Hermitian form $(~ ~,~ )$ is not positive
definite
on $F(n|m)$.

\n {\it Proof.} Consider a vector $|r,\0\r$. We recall, see (24),
that $0\le \phi <2\pi ,~~~\phi\ne \pi.$
For $\phi=0$ the RHS of (28) yields
$(|r,\0\r,|r,\0\r)=r!>0$.

If $\phi\ne 0$ (and $\phi\ne \pi$) (28) can be written as:
$$
(|r,\0\r,|r,\0\r)=\Big(cos(\phi/2)\Big)^{-r}\Big(sin(\phi)\Big)^{-r}
sin(\phi)sin(2\phi)sin(3\phi)\ldots sin(r\phi). \eqno(29)
$$
For $0<\phi <\pi $ both
$\Big(cos(\phi/2)\Big)^{-r}>0$ and $\Big(sin(\phi)\Big)^{-r}>0$.
Let $k$ be the smallest positive integer such that $0<k\phi<\pi$
and $\pi<(k+1)\phi<2\pi$.
Then $(|k,\0\r,|k,\0\r)>0$, but
$$
(|k+1,\0\r,|k+1,\0\r)=(|k,\0\r,|k,\0\r)
\Big(cos(\phi/2)\Big)^{-1}\Big(sin(\phi)\Big)^{-1}sin((k+1)\phi)<0.
\eqno(30)
$$

Since the replacement of $\phi$ with $-\phi$ does not change the RHS
of (30), $(|k+1,\0\r,|k+1,\0\r)<0$ also for $\pi<\phi<2\pi$.
This completes the proof.

\hfill $[]$

\bigskip
In view of the above proposition the Hermitian form (~,~) could be
positive definite eventually on $Cl_q(n|m)$-invariant subspaces
of $F(n|m) $ or in  factor spaces of them. 
For $q$ being a root of unity $F(n|m)$
contains such invariant subspaces. Indeed, set
$$
q=\exp ({i\pi l/k})\Longleftrightarrow \phi=\pi l/k,
\quad {\rm with}~~  l~~ {\rm and}~~ k ~~
{\rm being~relatively ~simple~integers},
\eqno(31)
$$
which (without loss of generality) we assume to be positive
integers and $l<k$. It is straightforward to verify that
$$
I_{l/k}(n|m)=span\{|r_1,\ldots,r_i,\ldots,r_{m+n}\r~|~
r_1\ge k,r_2\ge k,\ldots,r_m\ge k \} \eqno(32)
$$
is an (infinite-dimensional) invariant subspace.
The latter follows from (15c), namely
$$
c_j^-|r_1,\ldots,r_j,\ldots,r_{m+n}\r =
c_q [r_j] q^{r_1+\ldots+r_{j-1}}|\ldots,r_j-1,\ldots\r,
~~~j=1,\ldots,m \eqno(33)
$$
and the circumstance that $[r_j]=0$ for $r_j=k$
(and $q=\exp ({i\pi l/k})$).

Let $F_{l/k}(n|m)$ be the factor space of $F(n|m)$ with respect to
$I_{l/k}(n|m)$,
$$
F_{l/k}(n|m)=F(n|m)/I_{l/k}(n|m). \eqno(34).
$$
All vectors (considered as representatives of the equivalence
classes)
$$
|r_1,\ldots,r_{m+n}\r,\quad 0\le r_1,\ldots,r_m\le k-1,\quad 0\le
r_{m+1},
\ldots,r_{m+n}\le 1, \eqno(35)
$$
constitute a basis in $F_{l/k}(n|m)$.
For this reason we say that $F_{l/k}(n|m)$ can be considered as a subspace
of $F(n|m)$.

\vfill\eject
\n
{\bf Proposition 4.} The Fock space $F_{l/k}(n|m)$ is a
finite-dimensional
irreducible $Cl_q(n|m)-$module.

\n {\it Proof}.

(i) Let $j=1,\ldots,m$. Than (15c) yields
$$
c_j^-|r_1,\ldots,r_j,\ldots , r_{m+n}\r =\exp \big({i\pi{l\over
k}(r_1+\ldots+r_{j-1})}\big)
cos(\pi l/2k){sin(\pi{l\over k}r_j)\over{sin(\pi{l\over k})}}
|\ldots,r_j-1,\ldots\r. \eqno(36)
$$
Clearly $(-1)\exp \big({i\pi{l\over k}(r_1+\ldots+r_{j-1})}\big)\ne 0$. Since
$l<k\in \N$ ~ ($\N$ - all positive integers)
$cos(\pi l/2k)>0$ and $sin(\pi{l\over k})>0$. The
next observation is that  ${l\over k}r_j$ is not an
integer for $0<r_j<k$. Indeed, suppose that ${l\over k}r_j$ is an
integer
and let $s$ be the maximal integer common divisor of $k$ and $r_j$,
i.e.
$k=s.k'$ and $r_j=s.r'$.
Then ${l\over k}r_j={l.r'\over k'}$ and the
integers
$l.r'$ and $k'$ have no common divisors, i.e. ${l\over k}r_j$ is not
an
integer. Therefore $sin(\pi{l\over k}r_j)\ne 0$ and as a result the
coefficient
in front of $|\ldots,r_j-1,\ldots\r$ in the RHS of (36) is different
from zero
for any $r_j$, $0<r_j<k$.

(ii) If $m<j<m+n$ and $r_j=1$, then
$$
c_j^-|\ldots,r_j=1,\ldots \r =(-1)^{r_1+\ldots+r_{j-1}}
\exp \big({i\pi{l\over k}(r_1+\ldots+r_{j-1})}\big)
cos(\pi l/2k)
|\ldots,r_j=0,\ldots\r \eqno(37)
$$
and again the RHS of (37) is different from zero.

As a result
$$
(c_1^-)^{r_1}\ldots (c_{m+n}^-)^{r_{m+n}}|r_1,\ldots,r_{m+n}\r=
Const|0\r,\quad Const\ne 0. \eqno(38)
$$
From (38) and (15b) one easily concludes that $F_{l/k}(n|m)$ is an
irreducible $Cl_q(n|m)-$module. \hfill $[]$

So far we have defined a class of contravariant representations.
The problem that still  remains to be solved is
to select out of them those modules $F_{l/k}(n|m)$
for which the metric $(~,~)$ is positive definite, namely the modules,
which carry unitary representations.

\bigskip
\n {\bf Proposition 5.} The Fock space $F_{l/k}(n|m)$ is a
Hilbert space only if $l=1$ (and for any $1<k\in \N$), i.e. if $q$
is a primitive root of unity,
$$
q=\exp \big({i\pi/k}\big)\Longleftrightarrow \phi={\pi\over k}. \eqno(39)
$$
\n{\it Proof.}

\n (i) Assume $1<l<k.$

For $q=\exp ({i\pi l/k})$
and $1\le r < k$ one obtains from (28)
$$
(|r,\0\r,|r,\0\r)=\Big(cos({\pi\over 2}{l\over k})\Big)^{-r}
\Big(sin(\pi {l\over k})\Big)^{-r}
sin(\pi {l\over k}\, 1)sin(\pi {l\over k}\, 2)sin(\pi {l\over
k}\,3)\ldots
sin(\pi {l\over k}\, r). \eqno(40)
$$
Clearly,
$$
cos({\pi\over 2}{l\over k})>0,~~~ sin(\pi {l\over k})>0, ~~~{\rm
whereas}~~~
sin(\pi {l\over k}(k-1))<0.\eqno(41)
$$
Therefore there exists an integer $r_0,~~~1\le r_0 < k-1$ such that
$$
\pi {l\over k}\, 1< \pi,~~\pi {l\over k}\, 2< \pi,\ldots,
\pi {l\over k}\, r_0< \pi,~~{\rm whereas}~~
\pi {l\over k}\, (r_0+1) > \pi.\eqno(42)
$$
Then (40) yields $(|r_0+1\r,|r_0+1\r)<0$. This proves that the
Hermitian
form (~,~) is not a scalar product if $l\ne 1$ (if $q$ is not a
primitive
root of unity).

\bigskip
\n(ii) Let (39) holds. Then (28) yields
$$
(|r_1,\ldots , r_{m+n}\r ,|r_1,\ldots , r_{m+n} \r )=
cos({\pi\over 2k})^{-(r_1+\ldots+r_{m+n})}[r_1]!\ldots [r_j]!\ldots
[r_{m}]!.
\eqno(43)
$$
For each $j=1,\ldots,m$
$$
[r_j]!=\prod_{s=1}^{r_j}{sin({\pi\over k}\,s)\over{sin({\pi\over
k})}}>0
~~~{\rm since}~~~s\le r_j<k. \eqno(44)
$$
Since moreover $ cos({\pi/2k})>0$, the RHS of (43) is also
positive. This completes the proof.

\hfill $[]$

Clearly the dimension of the Fock space $F_{1/k}(n|m)$ is $k^m2^n$.
Define an orthonormal basis in $F_{1/k}(n|m)$
$$
\eqalignno{
&|r_1,\ldots , r_{m+n}) =\Big(c_q^{r_1+\ldots
+r_{m+n}} [r_1]![r_2]!\ldots [r_m]!\Big)^{-1/2}
|r_1,\ldots , r_{m+n}\r . &(45)\cr
}
$$
The transformation of this basis reads
($1<k\in \N$):
$$
\eqalignno{
&N_{j}|r_1,\ldots , r_{m+n}) =
r_j|r_1,\ldots , r_{m+n}),
 & (46a)\cr
&&\cr
&c_j^+|r_1,\ldots , r_{m+n}) =
\exp \Big({-i\pi (r_1+r_2+\ldots +r_{j-1})/k}\Big)
{\sqrt{2 sin(\pi (r_j+1)/k)sin (\pi /(2k))\over {sin^2(\pi/k)}}}&\cr
&\times |r_1, \ldots , r_{j-1}, r_j+1, r_{j+1}, \ldots, r_{m+n}),
~~\quad j=1,\ldots, m,
& (46b)\cr
& &\cr
&c_j^-|r_1,\ldots , r_{m+n})=
\exp \Big({i\pi (r_1+r_2+\ldots +r_{j-1})/k}\Big)
{\sqrt{2 sin(\pi r_j/k)sin (\pi /(2k))\over {sin^2(\pi/k)}}}&\cr
& \times |r_1, \ldots , r_{j-1}, r_j-1, r_{j+1}, \ldots, r_{m+n}),
~~\quad j=1,\ldots , m,
& (46c)\cr
& &\cr
&c_j^+|r_1,\ldots , r_{m+n}) =
(-1)^{r_1+\ldots +r_{j-1}}
(1-r_j)\exp \Big({-i\pi (r_1+\ldots +r_{j-1})/k}\Big)
{\sqrt{1\over {cos(\pi /(2k))}}} &\cr
&&\cr
&\times
|r_1, \ldots ,
r_{j-1}, r_j+1, r_{j+1},\ldots, r_{m+n} ), \quad j=m+1, \ldots ,m+n,
& (46d)\cr
&&\cr
&c_j^-|r_1,\ldots , r_{m+n}) =
(-1)^{r_1+\ldots +r_{j-1}}
r_j \exp( \Big({i\pi (r_1+\ldots +r_{j-1})/k}\Big)
{\sqrt{1\over {cos(\pi /(2k))}}}
 &\cr
&&\cr
&\times
|r_1, \ldots , r_{j-1}, r_j-1, r_{j+1},\ldots, r_{m+n} ),
\quad j=m+1, \ldots ,m+n.
& (46e)\cr
&&\cr
}
$$

The above relations yield that in the root of unity case the bosons 
are ``finite-dimensional'',
i.e. they are operators in finite-dimensional 
state spaces.

\bigskip\noindent
{\bf 3. A homomorphism of U$_{\bf q}$[osp(2n+1$|$2m)] onto 
Cl$_{\bf q}$(n$|$m). Applications to U$_{\bf q}$[osp(2n+1$|$2m)]
and U$_{\bf q}$[sl(m$|$n)] representations}

\bigskip\n
In what follows we 
observe 
that the deformed Bose and Fermi operators
(12) satisfy the defining relations (8) of the algebra  $\Uq$.
Therefore relations (46) give a root of unity representation also
of $\Uq$. The Cartan-Weyl elements of the subalgebra
$U_q[sl(m|n)]\subset \Uq $ are expressed in terms of the
$q$-anticommuting deformed Bose and Fermi operators. The
decomposition of  $\Uq$ Fock space $F_{1/k}(n|m)$ with respect to
$U_q[sl(m|n)]$ is considered, giving rise to root of unity
representations for the latter too.

\bigskip\n
{\bf Proposition 6.} The linear map
$\varphi:\Uq \rightarrow Cl_q(n|m)$,
defined on the generators of $\Uq$ as
$$
\varphi (a_i^{\pm})=c_i^{\pm}, ~~\varphi
(H_{i})=(-1)^{\l i \r}N_{i}-{1\over{2}}, ~~i=1,\ldots , m+n \eqno(47)
$$
and extended on all elements by associativity is a homomorphism
(in the sense of associative superalgebras) of $\Uq$ onto
$Cl_q(n|m)$.

\bigskip
The proof follows from the fact that the images $\varphi
(a_i^{\pm})$ and $\varphi(H_{i})$ of the generators
of $\Uq$ satisfy also the defining relations
(8). Moreover the generators of  $Cl_q(n|m)$ are among the images
of
$\varphi $: $c_i^{\pm}=\varphi (a_i^{\pm})$,
$N_i=\varphi((-1)^{\l i \r}(H_i+1/2))$.

\hfill $[]$

The relevance of Proposition 6 stems from the observation that
any $Cl_q(n|m)$ module $F$ and in particular the Fock modules
studied in the previous section can be turned into $\Uq$ modules
simply by setting $\varphi(a)\,x$ for any $a\in \Uq$ and $x\in
F$. The conclusion from (47) is that Eqs. (46) define an irreducible
$\Uq$ representation for any $q$ being a primitive root of unity.

\bigskip
Next we use the circumstance that the  quantum  superalgebra 
$U_q[sl(m|n)]$ is a subalgebra of $\Uq$. Therefore each
$\Uq$ module is (generally reducible) $U_q[sl(m|n)]$ module.
Using this we proceed to decompose  each $\Uq$ module $F_{1/k}(n|m)$
into a direct sum of irreducible $U_q[sl(m|n)]$ modules,
thus defining explicitly a class of root of unity representations
of the quantum  superalgebra $U_q[sl(m|n)]$. 

As a first step we recall the definition of
the quantum superalgebra $U_q[sl(m|n)]$. Let
$(\a_{ij})$ be the $(m+n-1) \times (m+n-1)$ Cartan matrix with
enties:
$$
\a_{ij}=(1+(-1)^{\t_{i,i+1}})\delta_{ij}-
(-1)^{\t_{i,i+1}}\delta_{i,j-1}-\delta_{i-1,j},\quad
i,j\in [1;n+m-1], \eqno(48)
$$
$$
\t_i=1-\l i \r =\cases {{\bar 0}, & if $\; i\leq m$,\cr
               {\bar 1}, & if $\; i> m$,\cr }; \quad
\t_{ij}=\t_i+\t_j. \eqno (49)
$$
{\bf Definition 4.} $U_q[sl(m|n)]$ is a Hopf algebra, which
is a topologically free $\Ch$ module, with Chevalley generators
$\hat{h}_i, \hat{e}_i, \hat{f}_i$, $ i=1,\ldots ,m+n-1,$ subject to
the
($deg(\hat{h}_i)={\bar 0},\quad
deg(\hat{e}_i)=deg(\hat{f}_i)=\t_{i,i+1}$)

\n
1. Cartan-Kac relations:
$$
\eqalignno{
& [\hat{h}_i,\hat{h}_j]=0,& (50a)\cr
& [\hat{h}_i,\hat{e}_j]=\a_{ij}\hat{e}_j,\quad
[\hat{h}_i,\hat{f}_j]
=-\a_{ij}\hat{f}_j,& (50b)\cr
& \[\hat{e}_i,\hat{f}_j\]=\delta _{ij}{\hat{k}_i-\bar{\hat{k}}_i
\over{q-\bar{q}}},\quad
  \hat{k}_i=q^{\hat{h}_i},\;\hat{k}_i^{-1}\equiv
{\bar{\hat{k}}}_i=q^{-\hat{h}_i}, & (50c)\cr
}
$$

\n
2. $\hat{e}$-Serre relations
$$
\eqalignno{
& [\hat{e}_i,\hat{e}_j]=0,\; if \; |i-j|\neq 1;\quad
{\hat{e}}^2_{m}=0,
& (51a) \cr
& [\hat{e}_i, [\hat{e}_{i}, {\hat{e}}_{i\pm 1}]_{\bar{q}}]_q=
  [\hat{e}_i, [\hat{e}_{i}, \hat{e}_{i\pm 1}]_{q}]_{\bar{q}}=0,
\quad i\neq m, & (51b)\cr
& \{ \hat{e}_{m},[[\hat{e}_{m-1},\hat{e}_{m}]_q,
\hat{e}_{m+1}]_{\bar{q}}\}=
   \{ \hat{e}_{m},[[\hat{e}_{m-1},\hat{e}_{m}]_{\bar{q}},
\hat{e}_{m+1}]_{q}\}=0, & (51c)\cr
}
$$

\n
2. $\hat{f}$-Serre relations, $\hfill (52)$

\n
obtained from the $\hat{e}$-Serre relations
by replacing everywhere
$\hat{e}_i$ with $\hat{f}_i$.

Setting (see (3)-(5))
$$
\eqalignno{
& h_i=\hat{h}_i, \;\;e_i=\hat{e}_i,\;\; i=1,\ldots , m,\;\; &\cr
& h_i=-\hat{h}_i,\;\; e_i=-\hat{e}_i, \; i=m+1, \ldots, m+n-1, & (53)
\cr
& f_i=\hat{f}_i, \;\; i=1, \ldots , m+n-1,&\cr
}
$$
one verifies that the defining relations of  $U_q[sl(m|n)]$
(50)-(52) are among the defining relations of  $\Uq$ (3)-(5).
Therefore $U_q[sl(m|n)]$ is a subalgebra of  $\Uq$ in a sense of
associative algebras.

A set of Cartan-Weyl elements of $U_q[sl(m|n)]$ has been
considered in [17], and consists of $m+n-1$ ``Cartan'' elements
$\tilde{H}_i=e_{11}-(-1)^{\theta_{i+1}}e_{i+1,i+1}=
\hat{h}_1+(-1)^{\t_2}\hat{h}_2+(-1)^{\t_3}\hat{h}_3+\ldots +
(-1)^{\t_i}\hat{h}_i$ and $(m+n)(m+n-1)$ root vectors $e_{ij}$,
$i\neq j=1,\ldots ,m+n$ ($\hat{e}_i=e_{i,i+1}$,
$\hat{f}_i=e_{i+1,i}$). $e_{ij}$ is positive if $i<j$ and
negative if $i>j$. Among the positive root vectors the normal
order is given by
$$
e_{ij}<e_{kl},~~
{\rm if }~~ i<k ~~{\rm  or }~~ i=k ~~{\rm and }~~ j<l; \eqno (54)
$$
for the negative root vectors $e_{ij}$ one takes the same
rule (54), and total order is fixed by choosing
$$
{\rm positive~ root~ vectors}<{\rm negative~ root~ vectors}<
\tilde{H}_{i}.
$$
A complete set of relations between the Cartan-Weyl elements is given
by:
$$
[ \tilde{H}_i, \tilde{H}_j]=0,  \eqno(55a)
$$
$$
[\tilde{H}_i , e_{jk}]=\big( \delta_{1j}-\delta_{1k}-(-1)^{\t_{i+1}}
(\delta_{i+1,j}-\delta_{i+1,k})\big) e_{jk}; \eqno(55b)
$$

\n
For any $e_{ij}>0$  and $e_{kl}<0$~:
$$
\eqalignno{
&\[e_{ij},e_{kl}\]={\Big (} (q-\bar{q}) \theta (j>k>i>l)(-1)^{\t_k}
e_{kj}e_{il}&\cr
& -\delta_{il}\theta (j>k)(-1)^{\t_{kl}}e_{kj}+
\delta_{jk}\theta (i>l) e_{il}{\Big )} \tL_{i-1}{\bar \tL}_{k-1}
&  \cr
& +\tL_{j-1}{\bar \tL}_{l-1}{\Big (} -(q-\bar{q})
\t (k>j>l>i) (-1)^{\t_j}
e_{il}e_{kj} &(55c)\cr
& -\delta_{il}\theta (k>j)(-1)^{\t_{ij}}e_{kj}+
\delta_{jk}\theta (l>i) e_{il}{\Big )} &\cr
& +{\delta_{il}\delta_{jk}\over {q-\bar{q}}}
\left( \tL_{j-1}^{(-1)^{\t_i}}{\bar \tL}_{i-1}^{(-1)^{\t_i}}-
{\bar \tL}_{j-1}^{(-1)^{\t_i}} \tL_{i-1}^{(-1)^{\t_i}}
\right); &\cr
}
$$
For $0<e_{ij}<e_{kl}$,
$$
\[e_{ij},e_{kl}\]_{q^{(-1)^{\t_j}\delta_{jl}-(-1)^{\t_j}\delta_{jk}
+(-1)^{\t_i}\delta_{ik}}}
=\delta_{jk}e_{il}+(q-\bar{q})(-1)^{\t_k}\t (l>j>k>i)
e_{kj}e_{il};\eqno(55d)
$$
For $0>e_{ij}>e_{kl}$,
$$
\[e_{ij},e_{kl}\]_{\bar{q}^{(-1)^{\t_j}\delta_{jl}-(-1)^{\t_j}\delta_{jk}
+(-1)^{\t_i}\delta_{ik}}}
=\delta_{jk}e_{il}-(q-\bar{q})(-1)^{\t_k}\t (i>k>j>l) e_{kj}e_{il},
\eqno(55e)
$$
$$
(e_{ij})^2=0 ~~ if ~~ \t_{ij}=1, \eqno(55f)
$$
where
$$
\theta (i_1>i_2>\ldots>i_k)=\cases{{1}, & if $i_1>i_2>\ldots>i_k$,
\cr
                             {0}, & otherwise \cr } ;
\quad \tL_i=q^{\tilde{H}_i}, ~~
\quad {\bar \tL}_{i}=\tL_i^{-1}.
$$
It is tedious but straightforward to check that the expressions:
$$
\eqalignno{
&e_{ij}={1\over 2}(-1)^{\l i\r}
\bar{q}^{(-1)^{\l j \r}N_j-1/2}
\[ c_i^-, c_j^+\],
~~i<j,
&  \cr
&e_{ij}={1\over 2}(-1)^{\l i\r} \[ c_i^-,
c_j^+\]
q^{(-1)^{\l i \r}N_i-1/2},
~~i>j,
& (56) \cr
& \tilde{H}_i=-N_1-(-1)^{\l i+1 \r}N_{i+1} &\cr
}
$$
satisfy the relations (55). Thus we have obtained a
realization
of the algebra $U_q[sl(m|n)]$ in terms of the $q$-anticommuting
deformed
Bose and Fermi operators (12) as a subalgebra of  $\Uq$.

The realization (56) allows one to consider any Fock module and
in particular the irreducible $\Uq$ Fock module $F_{1/k}(n|m)$ as
a $U_q[sl(m|n)]$ module and to decompose it into a direct sum of
irreducible $U_q[sl(m|n)]$ submodules.
From (56)
we have for the representatives $\rho (e_{ij})$, $i\neq j$, $\rho
(\tilde{H} _i)$ of all Cartan-Weyl generators of $U_q[sl(m|n)]$
in the Fock space $F_{1/k}(n|m)$:
$$
\eqalignno{
&\rho (e_{ij})=-(-1)^{\l i\r +\l j \r }
cos(\pi /(2k))
\exp\Big( -(-1)^{\l j \r}i\pi N_j/k\Big)
c_j^+ c_i^-, ~~i<j,
&  \cr
&\rho (e_{ij})=-cos(\pi /(2k))c_j^+ c_i^-
\exp\Big( (-1)^{\l i \r}i\pi N_i/k\Big ),
 ~~i>j,
& (57) \cr
& \rho (\tilde{H}_i)=-N_1-(-1)^{\l i+1 \r}N_{i+1}. &\cr
}
$$
The matrix elements follow from (46):
$$
\eqalignno{
& \rho (\tilde{H}_j)|r_1,\ldots , r_{m+n})=-\Big(
r_1+(-1)^{\l j+1\r }r_{j+1}\Big)
|r_1,\ldots , r_{m+n}), & (58a) \cr
&&\cr
& \rho (e_{jl})|r_1,\ldots , r_{m+n})=-\Big( 1-(1-\l l\r)r_l\Big)
(-1)^{(\l j\r-\l l \r)(r_1+\ldots +r_j)+(1-\l l \r)(r_{j+1}+\ldots
+r_{l-1})} \times
&\cr
&&\cr
& \exp\Big(-i\pi (r_j+\ldots +r_{l-1}+(-1)^{\l l\r }r_l-2\l l\r )/k
\Big)
{\sqrt{ sin(\pi r_j/k)sin (\pi(r_l+1) /k)\over {sin^2(\pi/k)}}}
|\ldots ,r_j-1,\ldots , r_l+1,\ldots ),&\cr
& \hskip 14cm ~~j<l,& (58b) \cr
&&\cr
& \rho (e_{jl})|r_1,\ldots , r_{m+n})=-\Big( 1-(1-\l l\r)r_l\Big)
(-1)^{(\l l\r-\l j \r)(r_1+\ldots +r_l)+(1-\l j \r)(r_{l+1}+\ldots
+r_{j-1})}\times
&\cr
&&\cr
&
\exp\Big( i\pi (r_l+\ldots +r_{j-1}+(-1)^{\l j\r }r_j)/k\Big)
{\sqrt{ sin(\pi r_j/k)sin (\pi(r_l+1) /k)\over {sin^2(\pi/k)}}}
|\ldots ,r_l+1,\ldots , r_j-1,\ldots ),&\cr
& \hskip 14cm ~~j>l.& (58c) \cr
}
$$
Obviously the subspace $F_r(m|n)$ spanned by all vectors (35)
with a fixed admissible value of $r_1+\ldots +r_{m+n}=r$ is an
$U_q[sl(m|n)]$ submodule. Since the matrix elements in front of
$|\ldots ,r_j-1,\ldots , r_l+1,\ldots )$ and $|\ldots
,r_l+1,\ldots , r_j-1,\ldots )$ in the RHS of (58) are different
from zero if these vectors belong to $F_{1/k}(n|m)$, the
$U_q[sl(m|n)]$ module $F_r(m|n)$ is irreducible. Hence the
irreducible $k^m2^n$ dimensional $\Uq$  Fock module
$F_{1/k}(n|m)$ splits into a direct sum of $mk-m+n+1$
irreducible  $U_q[sl(m|n)]$ modules  $F_r(m|n)$ characterized by
a number $r$ taking integer values from $0$ up to $m(k-1)+n$. For
different $k$ the modules $F_r(m|n)$ have different dimensions
and therefore they carry inequivalent irreducible primitive root
of unity representations of $U_q[sl(m|n)]$.

\bigskip\bigskip
\n
{\bf 4. Concluding remarks}

\bigskip\n
We introduced a construction for a deformation of the
Clifford superalgebra
$Cl(n|m)$.
As a guide for our construction 
we used the property  that
$Cl(n|m)$ is a homomorphic image of 
U[osp($2n+1|2m$)] and assumed 
that this property remains unchanged, i.e. that the
deformed superalgebra $Cl_q(n|m)$ must be  a homomorphic image of
$\Uq$. The resulting 
$Cl_q(n|m)$ is an associative superalgebra generated by $m$
deformed bosons and $n$ deformed fermions
with the somewhat unusual property that the bosons are
odd (i.e. fermionic) generators, whereas the fermions are even
(bosonic) generators.
A crucial role for determination of the
explicit form of the homomorphism $\varphi$ is played by the
realization of $\Uq$ via deformed Green generators.

The Fock space 
of $Cl_q(n|m)$ with the condition that
$$
\eqalignno{
&( c_i^+)^{\dag} =c_i^-, \quad ( N_{i})^{\dag} =N_{i},
\quad i=1,\ldots, m+n,
& \cr
}
$$
motivated from the requirement the physical 
observables to be hermitian operators
yields necessarily that  the deformed bosons are ``finite-dimensional'', 
i.e. they act
in finite-dimensional state spaces.

Every 
$F_{l/k}(n|m)$ is mapped by a homomorphism 
$\varphi $
into an irreducible $U_q[osp(2n+1|2m)]$ module, thus defining a
sequence of root of unity representations of this quantum superalgebra.
We decomposed each $F_{1/k}(n|m)$
into a direct sum of irreducible $U_q[sl(m|n)]$ modules and gave 
a class of root of 1 irreducible finite-dimensional 
representations
of $U_q[sl(m|n)]$. 
Finally, we derived  the action of all Cartan-Weyl elements $e_{ij}$ on 
the basis vectors. 

\vskip 30pt
\noindent
{\bf Acknowledgments.}

\s\n
TDP is  thankful to Prof. Randjbar-Daemi for the kind invitation
to visit the High Energy Section of the Abdus Salam International 
Centre for Theoretical Physics. NIS wishes to acknowledge the Alexander
von Humboldt Foundation for its support. This work was supported also
by the grant $\Phi-910$ of the Bulgarian Foundation for Scientific Research.

\vfill\eject\n
{\bf References}

\vskip 12pt
{\settabs \+  $^{11}\;\;$ & I. Patera, T. D. Palev, Theoretical
   interpretation of the experiments on the elastic \cr

\+ [1] & Khoroshkin S M  and Tolstoy V N 1991 Universal
         $R$-matrix for quantized (super)algebras \cr
\+     &   {\it Commun. Math. Phys.} {\bf 141} 599-617  \cr
\+ [2] & Palev T D 1998 A q-deformation of the parastatistics and an
alternative to the Chevalley description of  \cr
\+     &  $\Uq$ {\it Commun. Math. Phys.} {\bf 196} 429-443  \cr
\+ [3] & Palev T D 1982  Para-Bose and para-Fermi operators as generators
        of orthosymplectic Lie superalgebras \cr
\+    & {\it Journ. Math. Phys.} {\bf 23} 1100-1102 \cr
\+ [4] & Green H S 1953 A generalized method of field quantization
        {\it Phys. Rev.} {\bf 90} 270-273   \cr
\+ [5] & Kac V G and van der Leur J W 1988  Super boson-fermion
        correspondence of type B \cr
\+    & {\it Adv. Series Math.  Phys.} {\bf 7} 369-406 \cr
\+ [6] & Macfarlane A J 1989 On $q$-analogues of the quantum harmonic
      oscillator and the quantum group $SU(2)_q$ \cr
\+    & {\it J. Phys. A} {\bf 22} 4581-4588 \cr
\+ [7] & Biedenharn L C 1989 The quantum group $SU_q(2)$ and
   a $q$-analogue of the boson operators \cr
\+     & {\it J. Phys. A} {\bf 22} L873-L878 \cr
\+ [8] & Sun C P and Fu H C 1989 The $q$-deformed boson realisation
    of the quantum group $SU(n)_q$ and \cr
\+     & its representations {\it J. Phys. A} {\bf 22} L983-L986 \cr

\+ [9] & Hayashi T 1990 Q Analogs of Clifford and Weyl algebras: spinor and 
oscillator representations \cr
\+    & of quantum enveloping algebras  {\it  Commun.\ Math.\ Phys.} 
      {\bf 127}  129-144 \cr

\+ [10] & Pusz W and Woronowicz S L 1989  Twisted second quantization 
      {\it  Rep.\ Math.\ Phys.} 
      {\bf 27} 231-257 \cr

\+ [11] & Hadjiivanov L K, Paunov  R R and Todorov I T 1992
$U_q$ covariant oscillators and vertex operators  \cr        
\+     & {\it J.\ Math.\ Phys.} {\bf 33} 1379-1394  \cr

\+ [12]  & Jagannathan R, Sridhar R, Vasudevan R, Chaturvedi S,
          Krishnakumari M, Shanta P and \cr
\+      & Srinivasan V 1992 On the number operators of multimode 
systems of deformed 
oscillators covariant under \cr
\+       & quantum groups {\it  J.\ Phys.\ A} 
   {\bf 25} 6429-6454 \cr

\+ [13]  &Van der Jeugt J 1993 R-matrix formulation of deformed 
boson algebra  {\it  J.\ Phys.\ A~} 
     {\bf 26} L405-L411 \cr

\+ [14] & Floreanini R, Spiridonov V P and Vinet L 1991 q-Oscillator
        realizations of the quantum superalgebras \cr
\+     & $sl_q(m,n)$ and $osp_q(m,2n)$ {\it Commun. Math. Phys.}
         {\bf 137} 149-160  \cr
\+ [15] & Palev T D and Van der Jeugt J 1995 The quantum superalgebra
         $U_q[osp(1/2n)]$: deformed para-Bose  \cr
\+     & operators and root of
         unity representations  {\it J. Phys. A} {\bf   28} 2605-2616 \cr
\+ [16] & Kac V G and Raina A K 1987 Bombay lectures on highest weight 
representations of infinite dimensional \cr
\+     & Lie algebras  {\it Adv. Ser. Math. Phys.} {\bf 2} (World Scientific);
 Jacobson H P and Kac V G 1989 \cr
\+     & A New class of unitarizable highest weight representations
of infinite-dimensional Lie algebras \cr
\+     & {\it J. Funct. Anal.} {\bf 82} 69-90 \cr

\+ [17] & Palev T D, Stoilova N I and Van der Jeugt  J 2002 Jacobson
         generators of the quantum superalgebra \cr
\+    & $U_q[sl(n+1|m)]$ and Fock
        representations {\it Journ. Math. Phys.} {\bf 43} 1646-1663  \cr

\end